# Nonparametric "regression" when errors are positioned at end-points

PETER HALL[1,2] and INGRID VAN KEILEGOM[2]

[1]*Department of Mathematics and Statistics, University of Melbourne, Melbourne VIC 3010, Australia. E-mail: halpstat@ms.unimelb.edu.au*

[2]*Institut de Statistique, Université catholique de Louvain, Voie du Roman Pays 20, B-1348 Louvain-la-Neuve, Belgium. E-mail: ingrid.vankeilegom@uclouvain.be*

Increasing practical interest has been shown in regression problems where the errors, or disturbances, are centred in a way that reflects particular characteristics of the mechanism that generated the data. In economics this occurs in problems involving data on markets, productivity and auctions, where it can be natural to centre at an end-point of the error distribution rather than at the distribution's mean. Often these cases have an extreme-value character, and in that broader context, examples involving meteorological, record-value and production-frontier data have been discussed in the literature. We shall discuss nonparametric methods for estimating regression curves in these settings, showing that they have features that contrast so starkly with those in better understood problems that they lead to apparent contradictions. For example, merely by centring errors at their end-points rather than their means the problem can change from one with a familiar nonparametric character, where the optimal convergence rate is slower than $n^{-1/2}$, to one in the super-efficient class, where the optimal rate is faster than $n^{-1/2}$. Moreover, when the errors are centred in a non-standard way there is greater intrinsic interest in estimating characteristics of the error distribution, as well as of the regression mean itself. The paper will also address this aspect of the problem.

## 1. Introduction

The problem of estimating the end-point and tail shape of a distribution has a distinguished history, not least because it provides important examples of non-regular behaviour for various types of inference. See, for example, Harter and Moore (1965) and Smith (1985). The problem also has important practical motivations, arising in part from the prevalence of power-law distributions; see Zipf (1941, 1949). More recently, end-point and tail shape problems have been studied in regression settings; for example, in econometric models for auctions.

The importance of end-point estimation to auction models and the consequent fact that statistical inference in such models is non-regular were first noted by Paarsch (1992) and Donald and Paarsch (1993). The end-point problem arises there because the distribution of bid price generally depends on all the parameters of the model, for instance, on parameters that determine the costs of bidders. For particular examples of auction models, see Paarsch (1992) and Donald and Paarsch (2002).

Similar phenomena occur in truncated- or censored-regression models (e.g., Breen (1996); Long (1997)), market-structure analysis (e.g., Robinson and Chiang (1996)) and inference for production frontiers in econometrics (e.g., Aigner *et al.* (1977); Park and Simar (1994); Hall and Park (2004)). There is a strong association between these fields and those where extreme-value methods are used; for example, the successful bid at an auction is the extremum of all bids.

Although the term "regression" is commonly used in these settings, strictly speaking it is not correct. Since the error distribution is not centred at its expectation then the "regression mean" no longer admits its conventional definition as the average of the response variable given the value of the covariate, or explanatory, variable. This apparently minor distinction can have a major impact, and, for example, can lead to an intriguing paradox, as we shall show shortly.

In the context of auction models, Hirano and Porter (2003), Jofre-Bonet and Pesendorfer (2003) and Chernozhukov and Hong (2004) studied parametric approaches to inference about distribution end-points and jump heights. Campo *et al.* (2002) suggested a semi-parametric technique. Related statistical work tends to be in the setting of parametric regression; see, for example, Koenker *et al.* (1994), Smith (1994), Jurečková (2000), Portnoy and Jurečková (2000) and Knight (2001a). Knight (2001b) generalised several of the contributions of Smith (1994). In particular, in the context of curve estimation, he sketched the derivation of properties of estimators similar to those that we propose, although in cases where bias can be neglected and the shape parameter is fixed.

However, it is feasible to take a nonparametric view of this problem, permitting a greater degree of flexibility and generality. For example, Korostelev and Tsybakov (1993) treated a variety of boundary estimation problems from a nonparametric viewpoint. However, their work was generally in settings where information was available on both sides of the boundary and where data did not become relatively sparse as the boundary was approached. Therefore, the convergence rates derived by Korostelev and Tsybakov (1993) were faster than those that we give in the present paper. Chernozhukov (1998) addressed properties of nonparametric estimators alternative to ours, and derived upper bounds to convergence rates. Chernozhukov (2005) discussed a related problem, where a concise linear model, rather than a relatively unspecified smooth function, supplied the basis for inference. Among the contributions made here, over and above Chernozhukov's work, we give (for the somewhat different estimators proposed here) the structure of limiting distributions, provide a concise first-order description of the bias-variance tradeoff, discuss empirical choice of bandwidth, and give a detailed account of statistical issues, such as convergence rate paradoxes, for alternative, readily computable estimators.

The present paper suggests nonparametric methodology, and describes its properties, in the context of inference about end-point and tail shape functions in nonparametric regression. In this case the errors, or disturbances, in the nonparametric model are centred at their end-points, rather than at their means. The end-points may be assumed to take a convenient value such as zero. Thus, the problem of estimating the nonparametric regression mean becomes that of adaptively estimating the centring function.

Estimation of characteristics of the error distribution is sometimes also of practical interest. This problem can have several forms, depending on the extent of generality required. For example, if the error distribution has a jump discontinuity at its end-point then the height of the jump can be treated nonparametrically, or modelled parametrically, as a function of the explanatory variable. The end-point might be approached in a polynomial way, and then the exponent, or degree, may be one of the subjects of inference. This paper will address those issues, too.

The problem of nonparametric regression with end-point-centred errors also has significant theoretical motivation. In particular, depending on the way in which the end-point is approached, substantially faster convergence rates can be achieved than in conventional settings. For example, suppose we observe $Y_i = a(X_i) + \varepsilon_i$ for $1 \leq i \leq n$, where the errors $\varepsilon_i$ are independent and identically distributed with a distribution that has a jump discontinuity at one of its end-points and finite variance and $a$ denotes a twice-differentiable function. The estimator of $a$ given in this paper has root-mean-square convergence rate $n^{-2/3}$, which beats even the rate $n^{-1/2}$ for a parametric setting, let alone the rate $n^{-2/5}$ for standard nonparametric regression with twice-differentiable functions. We shall show that the rate $n^{-2/3}$ is minimax optimal.

However, it is well known that the rate $n^{-2/5}$ is also minimax optimal for estimating the same function. How can this be? This paradox can be resolved by noting that the two functions being estimated are not quite identical. They differ by a constant equal to the difference, $\delta$, between the mean and the end-point of the error distribution. The constant cannot be estimated at a faster rate than $n^{-2/5}$. However, this explanation is not without its own element of surprise, since it might be thought that estimation of $\delta$ would be a semi-parametric rather than a nonparametric problem; if we could observe the errors directly then we could estimate their end-point at rate $n^{-1}$ and their mean at rate $n^{-1/2}$, both expressed in root-mean-square terms.

## 2. Methodology

### 2.1. Model

Assume that data $(X_1, Y_1), \ldots, (X_n, Y_n)$ are generated by the model

$$Y_i = a(X_i) + \varepsilon_i, \qquad (2.1)$$

where $a$ denotes a smooth function, each $X_i$ is a $p$-vector and each $Y_i$ is a scalar. It is supposed that the distribution of the error, or disturbance, $\varepsilon_i$, conditional on $X_i = x$, has density $f(\cdot|x)$ with the property that $f(u|x) = 0$ for $u < 0$ and

$$f(u|x) = b(x)c(x)u^{c(x)-1} + \mathrm{O}(u^{c(x)+d-1}) \qquad \text{as } u \downarrow 0, \tag{2.2}$$

where $0 < d < \infty$. The quantities $b$ and $c$ are smooth, strictly positive functions from $\mathbb{R}^p$ to $\mathbb{R}$. We wish to estimate $a$, and sometimes also $b$ and $c$.

Taking the view that the locus of points $(x, a(x))$ represents a boundary, models similar to (2.1) and (2.2) have been treated before, although, with the exception of literature discussed in Section 1, generally only when $p = 1$. In the latter case, and in the context of statistics, contributions include those of Härdle *et al.* (1995), Hall *et al.* (1997, 1998), and especially Gijbels and Peng (2000). There is also a vast econometrics literature in the case $p = 1$, often including a shape constraint, such as convexity, in addition to smoothness of $a$. See, for example, Korostelev *et al.* (1995a, 1995b), Kneip *et al.* (1998) and Gijbels *et al.* (1999).

## 2.2. Nonparametric estimation of *a*

Let $h > 0$ denote a bandwidth. Given $x \in \mathbb{R}^p$, let $\mathcal{S}(x, h)$ be the set of pairs $(\alpha, \beta)$, where $\alpha$ is a scalar and $\beta$ is a $p$-vector, such that $Y_i \geq \alpha + \beta^{\mathrm{T}}(X_i - x)$ for all indices $i$ with $\|X_i - x\| \leq h$. Our initial estimator of $a(x)$ is

$$\tilde{a}(x) = \sup\{\alpha : (\alpha, \beta) \in \mathcal{S}(x, h)\}. \tag{2.3}$$

Should there be very few or no indices such that $\|X_i - x\| \leq h$, locally increase $h$ so that a moderate number of $X_i$'s lie in this range.

The one-sided nature of inference in this problem raises interesting issues connected with existence of the estimator and edge effects. To appreciate why, consider the case where the points $X_i$, for $1 \leq i \leq n$, all lie in a $p$-variate half-space defined by an infinite plane passing through $x$. Then there exists $\beta$ such that $\beta^{\mathrm{T}}(X_i - x) < 0$ for $1 \leq i \leq n$. Since the length of $\beta$ can be chosen arbitrarily large without altering the sign property, $\tilde{a}(x)$ as defined at (2.3) equals $+\infty$.

Let $\mathcal{R}$ denote the support of the common density, $g_X$, of the $X_i$'s, and write $\partial\mathcal{R}$ for the boundary of $\mathcal{R}$. If $g_X$ is continuous and positive in $\mathcal{R}$, and if $x$ is distant at least $sh$ (where $s > 0$) from $\partial\mathcal{R}$, then the probability that $\tilde{a}(x) = +\infty$ converges to zero exponentially fast, as a function of $n$, as the latter increases. See Section 5.1. However, if $x$ lies exactly on $\partial\mathcal{R}$, then, depending on the shape of the boundary, the probability can equal 1, even for finite $n$. Details are given in Section 3.1.

Arguably the simplest way of overcoming these difficulties is to set an upper bound, $B$, on the largest value that $a(x)$ can take and estimate $a(x)$ by averaging $\tilde{a}(u)$ over all values of $u$ for which $|x - u| \leq h_1$ and $|\tilde{a}(u)| \leq B$, where $h_1$ is another bandwidth. We shall discuss this approach in the next paragraph. Another method, more difficult to implement, is to distort the region of radius $h$ centred at $x$, within which $X_i$ must lie in order for $(X_i, Y_i)$ to be used to construct $\tilde{a}(x)$, so that the region includes values of $X_i$ that are further than $h$ from $x$ and appropriately complement the values of $X_i$ that are within $h$ of $x$.

One form that the averaging of $\tilde{a}(u)$ can take is based on local linear smoothing. There we choose $\hat{\alpha}_1 = \alpha_1 \in \mathbb{R}$ and $\beta_1 \in \mathbb{R}^p$ to minimise

$$\int \{\tilde{a}(x+h_1 u) - \alpha_1 - \beta_1^{\mathrm{T}} h_1 u\}^2 I\{|\tilde{a}(x+h_1 u)| \leq B\} K(u) \, \mathrm{d}u, \tag{2.4}$$

where $K$ is a bounded, spherically symmetric probability density supported on the $p$-variate unit sphere centred at the origin and $h_1$ is another bandwidth. Then we put $\hat{a}(x) = \hat{\alpha}_1$.

Alternatively, we may define

$$\check{a}(x) = \frac{\int_{\mathcal{R}(x)} \tilde{a}(x+h_1 u) I\{|\tilde{a}(x+h_1 u)| \leq B\} K(u) \, \mathrm{d}u}{\int_{\mathcal{R}(x)} I\{|\tilde{a}(x+h_1 u)| \leq B\} K(u) \, \mathrm{d}u}, \tag{2.5}$$

where $\mathcal{R}(x)$ denotes the set of points $u \in \mathbb{R}^p$ such that $x + h_1 u \in \mathcal{R}$, and $B$ is chosen sufficiently large to ensure that the denominator in (2.5) is non-vanishing. Both these approaches also overcome problems caused by discontinuities in the function $\tilde{a}$. While both address the issue of boundary effects, the estimator $\hat{a}$ suffers less from boundary bias than $\check{a}$. In both $\hat{a}$ and $\check{a}$ we may use a soft thresholding approach to inclusion of values of $u$ for which $|\tilde{a}(x+u)| \leq B$, rather than the hard thresholding suggested by (2.4) and (2.5). Since $\hat{a}$ and $\check{a}$ are based on local linear rather than local constant smoothing, they enjoy good performance near boundaries; our theoretical analysis in Section 3 will demonstrate this feature. General polynomial optimisation methods can also be employed to estimate $a$, although at the expense of greater computational labour.

Plug-in methods can be used to choose the bandwidth, $h$, empirically. However, motivation for that technique requires theory about large sample properties of $\tilde{a}$, and so discussion of empirical bandwidth selection is deferred to Sections 2.4 and 3.2.

The local linear estimator introduced in the first paragraph of this section can be viewed as based on a local, functional version of a linear programming algorithm. See Smith (1994) and Portnoy and Jurečková (2000) for related methodologies. The more general estimator, introduced in the paragraph above, requires polynomial programming for implementation.

If $c(x)$ lies in the interval $(0,2)$ then the rate of convergence of the estimator at (2.3) cannot be improved. If $c(x) \geq 2$ then the rate of convergence can be enhanced by using unboundedly many order statistics, where the number employed is a second smoothing parameter (in addition to the bandwidth $h$) and its optimal choice depends on knowing, or estimating, the main features of the remainder term in (2.2). See, for example, Hall (1982) and Smith (1985). However, when $c(x) \geq 2$ the data are very sparse in the neighbourhood of the boundary, and so inference about the remainder in (2.2) is especially difficult. Therefore, the empirical challenges posed by this approach usually outweigh any performance gains that might be achieved in practice, and so we shall not pursue such methods.

### 2.3. Nonparametric estimation of *b* and *c*

In principle, completely nonparametric methods may be used to estimate the functions $b$ and $c$, although in practice one would often take $c$ to be a constant, rather than a non-degenerate function of $x$.

When estimating $b$ and $c$ we need not use the numerical value of $h$ employed for $\tilde{a}$. However, in the brief account below we shall continue to use the notation $h$. Define the residuals $\tilde{\varepsilon}_i = Y_i - \tilde{a}(X_i)$, and let $\mathcal{T}(x,h)$ denote the set of $\tilde{\varepsilon}_i$'s for which $\tilde{\varepsilon}_i > 0$ and $\|X_i - x\| \leq h$. Put $N_1 = \#\mathcal{T}(x,h)$, and rank the elements of $\mathcal{T}(x,h)$ as $0 < \hat{\varepsilon}_{(1)}(x,h) \leq \cdots \leq \hat{\varepsilon}_{(N_1)}(x,h)$. Put

$$\hat{c}(x) = \left\{\log \hat{\varepsilon}_{(r+1)}(x,h) - \frac{1}{r}\sum_{i=1}^{r} \log \hat{\varepsilon}_{(i)}(x,h)\right\}^{-1}, \qquad \hat{b}(x) = (r/N_1)\{\hat{\varepsilon}_{(r+1)}(x,h)\}^{-\hat{c}(x)},$$

where $r$, another smoothing parameter, denotes a threshold. Optimal choice of bandwidth for estimating $b$ and $c$ is a highly complex matter, and will not be treated here. The estimators $\hat{b}$ and $\hat{c}$ can be thought of as local function versions of conditional maximum likelihood estimators suggested by Hill (1975).

### 2.4. Outline of theoretical properties

We shall show in Section 3 that, when constructing the local linear estimator $\tilde{a}$ and its smoothed versions $\hat{a}$ and $\check{a}$, it is generally optimal to choose $h \sim const.n^{-1/(p+2c)}$. In this case the estimators have root-mean-square convergence rate $n^{-2/(p+2c)}$, when applied to cases where $a$ has two derivatives. For very general choices of the error distribution, this rate is optimal when $0 < c < 2$. Even if the functions $b$ and $c \in (0,2)$ take known, constant values, and we know the error distribution exactly (e.g., that it is gamma or Weibull), the rate $n^{-2/(p+2c)}$ cannot be improved upon.

However, when $c \geq 2$, and we have sufficient information about the error distribution, the convergence rate of estimators of $a$ can be improved by using other approaches. For instance, if $b$ and $c$ are constant, and if the error density $f$ is known, then an estimator of $a$ that is based on maximising a "local" version of log-likelihood can produce an estimator that converges to $a$ at rate $n^{-2/(p+4)}$, rather than $n^{-2/(p+2c)}$, when $p > 2$ and $a$ has two derivatives.

The problem is more awkward when the error distribution is not known. There, the convergence rate $n^{-2/(p+2c)}$ can be close to optimal. In particular, if we know only that the errors have a common density $f$, with $f(u) = bcu^{c-1} + \mathrm{O}(u^{c+d})$ as $u \downarrow 0$, where $b, c > 0$ are fixed constants, then the minimax optimal convergence rate of estimators of $a$ is $n^{-2/(p+2c)-\delta(d)}$, where $\delta(d) > 0$ converges to zero as $d \downarrow 0$.

## 3. Theoretical properties

### 3.1. Convergence rates of estimators of *a*

Assume that data $(X_i, Y_i)$ are generated by the model at (2.1), where

> the pairs $(X_1, \varepsilon_1), (X_2, \varepsilon_2), \ldots$ are independent and identically distributed as $(X, \varepsilon)$; the density of $X$ is supported in a compact region, $\mathcal{R} \subseteq \mathbb{R}^p$, and is continuous and non-zero there; $P(\varepsilon > 0) = 1$; the distribution of $\varepsilon$, conditional on $X = x$, is absolutely continuous with a density, $f(\cdot|x)$, which satisfies (2.2) and, in the notation there, $d > 0$ is fixed, $b$ and $c$ are Hölder-continuous functions satisfying $C_1 \leq b(x), c(x) \leq C_2$ for all $x \in \mathcal{R}$, $C_1$ and $C_2$ are constants satisfying $0 < C_1 < C_2 < \infty$ and the remainder in (2.2) is of the order stated there, uniformly in $x \in \mathcal{R}$; and $\sup_x E(\varepsilon^{2+\eta}|X = x) < \infty$ for some $\eta > 0$. (3.1)

Recall from Section 2 that the one-sided nature of the inference problem means that the estimator $\tilde{a}$ will often tend not to be defined at the boundary. However, $\tilde{a}$ may be well-defined very close to the boundary. To elucidate this behaviour we shall consider two types of $x$, described in (3.2) below. By way of notation, given a point $x_0$ in the boundary $\partial \mathcal{R}$ of $\mathcal{R}$, let $v(x_0)$ denote the inward-pointing normal to the tangent plane at $x_0$, which is well defined if $\partial \mathcal{R}$ has a continuously turning tangent in a neighbourhood of $x_0$. Then we ask that:

> Either $x$ is fixed as an interior point of $\mathcal{R}$, or $x = x(n)$ is within order $h$ of $\partial \mathcal{R}$, in the following sense: Suppose that for some $x_0 \in \partial \mathcal{R}$, $\partial \mathcal{R}$ is of codimension 1 and has a continuously turning tangent in a neighbourhood of $x_0$, and that for some sufficiently small $\delta > 0$, $x_1 + v(x_1)t \in \mathcal{R}$ for all $x_1 \in \partial \mathcal{R}$ with $\|x_0 - x_1\| \leq \delta$ and all $0 \leq t \leq \delta$. Define $x = x(n) = x_0 + v(x_0)sh$, where $x_0 \in \partial \mathcal{R}$, $s > 0$ and $x_0$ and $s$ are held fixed. (3.2)

If $x \in \mathcal{R}$ is an interior point, or if $x = x(n) = x_0 + v(x_0)sh$ where $x_0 \in \partial \mathcal{R}$ and $s \geq 1$, let $\mathcal{U}(x)$ denote the closed, $p$-variate sphere of unit radius centred at $x$. If $x = x(n) = x_0 + v(x_0)sh$ with $x_0 \in \partial \mathcal{R}$ and $0 < s < 1$, take $\mathcal{U}(x)$ to be the larger of the two parts of the just-mentioned sphere that are obtained by cutting it by the plane that is perpendicularly distant $s$ from the origin and has its normal in the direction $v(x)$, pointing towards the centre to the sphere.

Let $\dot{a}$ and $\ddot{a}$ denote the $p$-vector of first derivatives and $p \times p$ matrix of second derivatives of the function $a$ and suppose that

> the function $a$ has two continuous derivatives in $\mathcal{R}$, and if $x = x_0 + v(x_0)sh$ then $\partial \mathcal{R}$ has a continuously turning tangent plane at $x_0$. (3.3)

Assume, too, that

for some $0 < \eta < 1/(2p)$ and all sufficiently large $n$, $n^{\eta-(1/p)} < h < n^{-\eta}$. (3.4)

Given $x \in \mathcal{R}$, let $E_1, E_2, \ldots$ denote independent, exponentially distributed random variables (all with unit mean), write $\gamma$ for Euler's constant and define

$$Z_j(x) = \exp\left[-c(x)^{-1}\left\{\sum_{i=j}^{\infty}(E_i - 1)i^{-1} + \gamma - \sum_{i=1}^{j-1} i^{-1}\right\}\right], \quad j \geq 1. \quad (3.5)$$

Given $x \in \mathcal{R}$, let $U_1(x), U_2(x), \ldots$ be independent and identically distributed random $p$-vectors, independent of the $Z_j(x)$'s and uniformly distributed on $\mathcal{U}(x)$. For $c_1, c_2 \geq 0$, define

$$Q_1(c_1, c_2|x) = \sup_{\beta \in \mathbb{R}^p} \inf_{1 \leq i < \infty} [c_1\{\beta^{\mathrm{T}} U_i(x) + \tfrac{1}{2} U_i(x)^{\mathrm{T}} \ddot{a}(x) U_i(x)\} + c_2 b(x)^{-1/c(x)} Z_i(x)].$$

Note that $Q(1, 0|x)$ is a constant. In Theorem 1, below, this degenerate distribution is a limit in the case where $\tilde{a}(x) - a(x)$ is asymptotically dominated by bias.

In the statement of Theorem 1 we let $x_1$ denote $x$ if $x$ is an interior point of $\mathcal{R}$, and $x_1 = x_0$ if $x = x(n) = x_0 + v(x_0)sh$. Let $w(p)$ be the content of the $p$-variate unit sphere (thus, $w(1) = 2$, $w(2) = \pi$), let $g_X(x)$ represent the value of the density of the distribution of $X$ at $x$ and put $w_x = w(p) g_X(x_1)$. (To simplify notation we suppress the role of $x_1$ here.) We use a simpler rule than that in Section 2.2 to take care of cases where $\tilde{a}(x)$ is infinite. However, the last sentence in the theorem remains true if we define $\tilde{a}(x)$ to equal zero whenever $|\tilde{a}(x)| > B$, provided $B > |a(x)|$.

**Theorem 1.** *Assume (3.1)–(3.4). (a) If $(w_x n h^p)^{1/c(x_1)} h^2 \to \rho$, where $\rho \in [0, \infty)$, then $(w_x n h^p)^{1/c(x_1)}\{\tilde{a}(x) - a(x)\} \to Q_1(\rho, 1|x_1)$ in distribution. (b) If $(nh^p)^{1/c(x_1)} \times h^2 \to \infty$ then $h^{-2}\{\tilde{a}(x) - a(x)\} \to Q_1(1, 0|x_1)$ in distribution. Furthermore, if we take the precaution of defining $\tilde{a}(x)$ to equal an arbitrary but fixed constant in cases where it would otherwise be infinite, then second moments converge to those of the limiting distributions.*

Proofs of Theorems 1–3 will be given in Section 5. It is crucial, in condition (3.2), that we take $s > 0$ rather than $s \geq 0$. If $s = 0$ then $x$ lies right on the boundary of $\mathcal{R}$, and in such cases the theorem is false. For example, if $\mathcal{R}$ is a convex region with a smooth boundary, such as a sphere, then with probability 1, $\tilde{a}(x_0) = \infty$ for all $x_0 \in \partial \mathcal{R}$. However, it follows from the theorem that for points $x$ that are arbitrarily close to $\partial \mathcal{R}$, on the scale of the bandwidth, without being right on the boundary, the probability that $\tilde{a}(x)$ is finite converges to 1, and in fact the estimator $\tilde{a}(x)$ attains optimal convergence rates. The main impact of the boundary is to reduce the number of design points used to construct the estimator. This tends to inflate estimator variance, although only by a constant factor (determined by the geometry of the boundary in the vicinity of $x$). The difficulty can be alleviated by increasing the bandwidth in such places.

Asymptotic properties of $\hat{a}$ and $\check{a}$ are similar, except that the limiting distribution of $\hat{a}$ is more tedious to define. Therefore we shall confine ourselves to $\check{a}$. To further abbreviate our treatment we shall restrict attention to the case where

$$x \text{ is an interior point of } \mathcal{R}, \; h_1 = th \text{ for a fixed constant } t > 0, \text{ and} \\ (w_x n h^p)^{1/c(x)} h^2 \to \rho \in [0, \infty). \tag{3.6}$$

Let $Z_1, Z_2, \ldots$ be as at (3.5); for simplicity we drop the argument $x$. Re-define $U_1, U_2, \ldots$ to be independent of one another and of the $Z_j$'s and uniformly distributed in the $p$-variate sphere of radius $t+1$ centred at $x$. Given a $p$-vector $u$ with $\|u\| \leq t$, let $(S_1(u), T_1(u)), (S_2(u), T_2(u)), \ldots$ denote the values $(U_{i_1(u)}, Z_{i_1(u)}), (U_{i_2(u)}, Z_{i_2(u)}), \ldots$ of $(U_i, Z_i) = (U_i, Z_i(x))$ for which $\|U_i - u\| \leq h$, ordered such that $Z_{i_1(u)} < Z_{i_2(u)} < \cdots$. With $\kappa = p^{-1} \int \|u\|^2 K(u) \, du$, $\nabla^2$ denoting the Laplacian operator and $\rho \geq 0$ as in (3.6), define

$$Q_2(u|x) = \sup_{\beta \in \mathbb{R}^p} \inf_{1 \leq i < \infty} [\rho \{ \beta^{\mathrm{T}} S_i(u) + \tfrac{1}{2} S_i(u)^{\mathrm{T}} \ddot{a}(x) S_i(u) \}$$

$$+ \{(t+1)^p b(x)\}^{-1/c(x)} T_i(u)],$$

$$Q_3(x) = \frac{1}{2} \rho t^2 \kappa (\nabla^2 a)(x) + \int Q_2(u|x) K(u) \, du.$$

Under conditions (3.1)–(3.6), and taking $B > |a(x)|$ in (2.5), it can be shown that with probability $1 - \mathrm{O}(n^{-C})$ for all $C > 0$, the estimator $\check{a}(x)$ at (2.5) satisfies

$$\check{a}(x) = \int \tilde{a}(x + h_1 u) K(u) \, du. \tag{3.7}$$

Theorem 2 applies with equal validity to the estimators at (2.5) and (3.7). Nevertheless, the estimator on the right-hand side of (3.7) does not enjoy the boundedness property that partly motivated $\check{a}$ at (2.5).

**Theorem 2.** *Assume (3.1)–(3.6), and that the kernel $K$ used to define $\check{a}(x)$ is a bounded, spherically symmetric probability density supported on the unit sphere centred at the origin. Then $(w_x n h^p)^{1/c(x)} \{\check{a}(x) - a(x)\} \to Q_3(x)$ in distribution. Furthermore, if in the integrand at (3.7) we take the precaution of defining $\check{a}(x + h_1 u)$ to equal an arbitrary but fixed constant in cases where it would otherwise be infinite, then the second moment converges to that of the limiting distribution.*

### 3.2. Choice of bandwidth

Theorems 1 and 2 imply that, except in pathological cases where $\ddot{a}(x) = 0$, the optimal convergence rate of $\tilde{a}(x)$ and $\check{a}(x)$ to $a(x)$ is achieved by choosing the bandwidth $h$ so that $(nh^p)^{-1/c(x)}$ and $h^2$ are of the same size and, in particular, $h \sim const.n^{-1/\{p+2c(x)\}}$.

If $x$ does not lie on the boundary of $\mathcal{R}$, and if $(w_x n h^p)^{1/c(x)} h^2 \to \rho \in [0, \infty)$, then the asymptotic mean squared error of $\tilde{a}(x)$ is given by

$$\tau(\rho|x) = E\Big\{\sup_{\beta \in \mathbb{R}^p} \inf_{1 \leq i < \infty} [\rho\{\beta^{\mathrm{T}} U_i + \tfrac{1}{2} U_i^{\mathrm{T}} \ddot{a}(x) U_i\} + b(x)^{-1/c(x)} Z_i(x)]\Big\}^2, \qquad (3.8)$$

where $U_1, U_2, \ldots$ are uniformly distributed on the unit sphere centred at $x$, $Z_1(x), Z_2(x), \ldots$ are defined at (3.5) and the $U_i$'s and $Z_i(x)$'s are completely independent. Therefore, if $h = w_x^{-1/\{p+2c(x)\}} \rho^{1/\{2+p/c(x)\}} n^{-1/\{p+2c(x)\}}$ then $\rho$ should ideally be chosen as

$$\rho_0(x) = \arg\min_\rho \tau(\rho|x). \qquad (3.9)$$

One technique for estimating $\ddot{a}(x)$ is to twice numerically differentiate a heavily smoothed version of $\tilde{a}$. A simpler approach, if we may make the assumption (A), say, that, for each $i$, the distribution of $\varepsilon_i$ does not depend on $X_i$, is to pass a traditional smoother through the data $(X_i, Y_i)$ estimates the value of $\mu(x) = E(Y|X = x)$. Under (A), this quantity differs from $a$ only by a constant, and so $\ddot{a} = \ddot{\mu}$. The latter function can be estimated using conventional cubic smoothing. This approach is attractive even if the distribution of $\varepsilon_i$ depends to some extent on $X_i$, since it gives a working empirical approximation to $\ddot{a}$.

Methods for estimating $b(x)$ and $c(x)$ were discussed in Section 2.3. Substituting these estimators for the true values of $\ddot{a}(x)$, $b(x)$, $c(x)$ and $\rho(h, x)$ in (3.8), we may compute an estimator $\hat{\tau}(\rho|x)$ of $\tau(\rho|x)$ using a Monte Carlo simulation, which leads to an estimator $\hat{\rho}_0(x)$ of $\rho_0(x)$ at (3.9). The density of $X$ at $x$, that is, $g_X(x)$, can be estimated more conventionally, and thus an estimator $\widehat{w}_x$ of $w_x = w(p) g_X(x)$ can be constructed. An empirical bandwidth selector is then given by

$$h(x) = \widehat{w}_x^{-1/\{p+2\hat{c}(x)\}} \hat{\rho}_0(x)^{1/\{2+p/\hat{c}(x)\}} n^{-1/\{p+2\hat{c}(x)\}}. \qquad (3.10)$$

In many circumstances it is feasible to take $c(x)$ to be a constant, not depending on $x$. Then a global approach to bandwidth choice is possible, as follows: We shall proceed as though the density $g_X$ is constant; if it is not, using its average value rather than attempting to accommodate its variation greatly simplifies matters. Thus, we take $\widehat{w}$ to be an estimator of the average value of $w_x$. The mean integrated squared error of $\tilde{a}(x)$ is asymptotic to $\tau(\rho) = \int_\mathcal{R} \tau(\rho|x) \, dx$, of which an estimator is $\hat{\tau}(\rho) = \int_\mathcal{R} \hat{\tau}(\rho|x) \, dx$, leading to an estimator $\hat{\rho}_0 = \arg\min_\rho \hat{\tau}(\rho)$ of $\rho_0 = \arg\min_\rho \tau(\rho)$. A global bandwidth for constructing $\tilde{a}$ is thus $h = \widehat{w}^{-1/(p+2\hat{c})} \hat{\rho}_0^{1/(2+p/\hat{c})} n^{-1/(p+2\hat{c})}$.

### 3.3. Optimality

We shall show in this section that the convergence rates implied by Theorems 1 and 2, and also lower bounds of the same orders, are available uniformly over classes $\mathcal{A}$ of functions $a$ with two bounded derivatives. The possibility that either the proportionality constant, $b$, or the exponent, $c$, varies with the design variable, $X_i$, is not relevant to discussion of the lower bound, and for this reason, for simplicity, and since our lower bound results are stronger if we narrow the class of error distributions for which worst-case performance is achieved, we shall take the distribution of $\varepsilon = \varepsilon_i$ to be a single, specific one, say the gamma:

$$f(u) = f_i(u) = \frac{1}{\Gamma(c)} u^{c-1} e^{-u}, \qquad \text{where } c > 0 \text{ is fixed.} \tag{3.11}$$

In the lower bound calculations, $c > 0$ will be assumed known.

Likewise, we shall treat just one distribution of $X = X_i$ and one region $\mathcal{R}$. In particular, writing $\mathcal{V}(x, r)$ for the closed sphere centred at $x$ and of radius $r > 0$, we shall assume that

$$\mathcal{R} = \mathcal{V}(x_0, 1) \quad \text{and} \quad X \text{ is uniformly distributed on } \mathcal{R}. \tag{3.12}$$

Given $C > 0$, let $\mathcal{A} = \mathcal{A}(C)$ denote the class of functions $a$ for which first and second derivatives exist and are bounded absolutely by $C$, let $\bar{\mathcal{A}}$ denote the class of bounded functions $\bar{a}$ of the data $(X_1, Y_1), \ldots, (X_n, Y_n)$ (the latter generated as at (2.1)) and let $\mathcal{R}_h$ be the set of all points in $\mathcal{R}$ that are distant at least $h$ from $\partial \mathcal{R}$.

**Theorem 3.** *Assume (3.11) and (3.12) and, when constructing $\tilde{a}(x)$, let $h = \text{const.} \times n^{-1/(p+2c)}$, except that we take $\tilde{a}(x)$ equal to an arbitrary but fixed constant in cases where it would otherwise be infinite. Then,*

$$\sup_{x \in \mathcal{R}_h} \sup_{a \in \mathcal{A}} E\{\tilde{a}(x) - a(x)\}^2 = \mathrm{O}(n^{-2/(p+2c)}) \tag{3.13}$$

*as $n \to \infty$. Furthermore, if $0 < c < 2$,*

$$\liminf_{n \to \infty} n^{2/(p+2c)} \inf_{\bar{a} \in \bar{\mathcal{A}}} \sup_{a \in \mathcal{A}} E\{\bar{a}(x) - a(x)\}^2 > 0 \qquad \text{for each } x \in \mathcal{R} \setminus \partial \mathcal{R}, \tag{3.14}$$

$$\liminf_{n \to \infty} n^{2/(p+2c)} \inf_{\bar{a} \in \bar{\mathcal{A}}} \sup_{a \in \mathcal{A}} \int_{\mathcal{R}_h} E\{\bar{a}(x) - a(x)\}^2 \, \mathrm{d}x > 0. \tag{3.15}$$

Together, (3.13)–(3.15) imply that the estimator $\tilde{a}$ achieves the minimax optimal rate, $n^{-2/(p+2)}$, uniformly over all functions $a \in \mathcal{A}$, and that the optimality can be expressed in either a local or a global sense. Similarly, it may be proved that if $\mathcal{A}_q$ is taken to be the class of functions $a$ with $q + 1$ (rather than 2) bounded derivatives, then the $q$th degree local polynomial approach discussed in Section 2.2 achieves the minimax optimal convergence rate of $n^{-2q/(p+2cq)}$ uniformly over functions in $\mathcal{A}_q$. The upper bound (3.13) continues to hold if the class $\mathcal{A}$ is increased to include a range of distributions of $\varepsilon$ for which the lower tail of the distribution function decreases like $u^c$ as $u \downarrow 0$, and a range of distributions of design points for which the density is bounded away from zero on $\mathcal{R}$.

# 4. Numerical properties

## 4.1. Simulations

Consider independent and identically distributed data $(X_i, Y_i)$ $(1 \leq i \leq n)$ satisfying the model $Y_i = a(X_i) + \varepsilon_i$ given in (2.1). The covariate $X_i$ has a uniform distribution on the interval $[0, 1]$. We consider three models for $a(x)$ $(0 \leq x \leq 1)$:

$$
\begin{aligned}
&\text{Model 1}: a(x) = 10(x - a_0)^3, &&a_0 = 0.25, 0.5,\\
&\text{Model 2}: a(x) = \exp(-a_0 x^2), &&a_0 = 1, 2,\\
&\text{Model 3}: a(x) = a_0 \cos(\pi x), &&a_0 = 0.25, 0.5.
\end{aligned}
\quad (4.1)
$$

Figure 1 shows the graphs of these six frontier functions. The error $\varepsilon_i$ is taken from a Gamma distribution:

$$ f(u|x) = \frac{1}{s(x)^c \Gamma(c)} u^{c-1} \exp\{-u/s(x)\} $$

$(u \geq 0)$, where $c > 0$ and $s(x) = 1 + 2x$. Note that this density is of the general type (2.2), with $b(x) = \{cs(x)^c \Gamma(c)\}^{-1}$.

We carry out two simulation studies. In the first study, we investigate the performance of the estimator $\hat{a}(x)$ and of its data-driven bandwidth selector described in Section 3.2. The simulations are executed based on 100 arbitrary samples of size $n = 200$ and $n = 400$. For each sample we estimate $a(x)$ at $x = 0.5$. The scaling parameter $c$ is chosen to be 0.5, 1 or 1.5. These three values of $c$ are such that, as $u \downarrow 0$, $f(u|x) \to \infty$, $f(u|x) \to s(x)^{-1}$ and $f(u|x) \to 0$, respectively. We use local linear smoothing to obtain both $\tilde{a}(x)$ and $\hat{a}(x)$. The bandwidth $h$ is calculated from formula (3.10) and we have taken $h_1 = h$. To estimate $\ddot{a}(x)$ we work (as explained in Section 3.2) under the working model that the distribution of $\varepsilon_i$ does not depend on $X_i$, in which case $\ddot{a}(x)$ equals the second derivative of the regression function $E(Y|X=x)$. This second derivative is estimated using local cubic smoothing, with bandwidth 0.25. The functions $b(x)$ and $c(x) \equiv c$ are estimated employing the procedure explained in Section 2.3, where $r$ equals the smallest integer larger than $0.90 N_1$ and the bandwidth for estimating $b(x)$ and $c(x)$ is chosen as 0.25. The kernel used throughout is the biquadratic kernel, $K(u) = (15/16)(1-u^2)^2 I(|u| \leq 1)$.

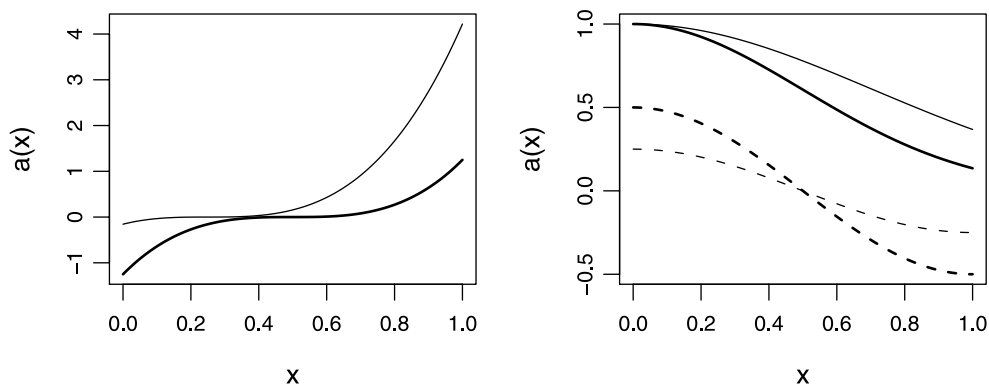

**Figure 1.** Graphs of the functions $a(x)$ given in (4.1): the left figure shows $a(x)$ for Model 1 ($a_0 = 0.25$ (thin curve) and $a_0 = 0.50$ (thick curve)), the right figure shows $a(x)$ for Model 2 ($a_0 = 1$ (thin solid curve) and $a_0 = 2$ (thick solid curve)) and Model 3 ($a_0 = 0.25$ (thin dashed curve) and $a_0 = 0.50$ (thick dashed curve)).

Tables 1 and 2 show the estimated bias, variance and mean square error (MSE) of $\hat{a}(x)$ at $x=0.5$ for each of the considered models as well as the average value of the bandwidth $h$ over the 100 simulation runs obtained using a Monte Carlo simulation of formula (3.10). Note that the functions $a(x)$ considered in this simulation study are neither convex nor concave. In fact, our method imposes neither condition in contradistinction to, for instance, the DEA (data envelopment analysis) estimator, which requires the function $a(x)$ to be convex.

The tables show that the MSE increases when $c$ increases, which is to be expected since the higher the value of $c$, the smaller the density $f(\cdot|x)$ of the error close to the frontier, and so the harder the estimation of the frontier. These findings also agree with the theoretical results of Section 3. This sparsity of data close to the frontier affects especially the bias of the estimator, since it is clear that the estimator $\tilde{a}(x)$ tends to overestimate $a(x)$ whenever there are few observations near the boundary. It also affects, as would be expected, the performance of the empirical bandwidth selector. The higher the value of $c$, the harder it becomes to estimate in an accurate way the optimal bandwidth. Finally, comparing Tables 1 and 2 we see that both the bias and the variance decrease as the sample size increases.

**Table 1.** Monte Carlo simulations for $n = 200$, with optimal bandwidth given by (3.10)

| Model | $a_0$ | $c$ | Mean($h$) | 10 Bias | 100 Var | 100 MSE |
|---|---|---|---|---|---|---|
| 1 | 0.25 | 0.5 | 0.045 | 0.156 | 0.073 | 0.098 |
|   |      | 1   | 0.071 | 1.408 | 1.244 | 3.226 |
|   |      | 1.5 | 0.096 | 3.177 | 1.905 | 11.995 |
|   | 0.5  | 0.5 | 0.067 | −0.169 | 0.049 | 0.078 |
|   |      | 1   | 0.094 | 0.798 | 0.707 | 1.344 |
|   |      | 1.5 | 0.120 | 2.389 | 1.832 | 7.540 |
| 2 | 1    | 0.5 | 0.064 | −0.010 | 0.023 | 0.023 |
|   |      | 1   | 0.087 | 0.899 | 0.812 | 1.619 |
|   |      | 1.5 | 0.119 | 2.406 | 1.954 | 7.741 |
|   | 2    | 0.5 | 0.051 | 0.079 | 0.033 | 0.039 |
|   |      | 1   | 0.079 | 1.072 | 1.032 | 2.181 |
|   |      | 1.5 | 0.111 | 2.551 | 1.795 | 8.300 |
| 3 | 0.25 | 0.5 | 0.062 | 0.009 | 0.022 | 0.022 |
|   |      | 1   | 0.086 | 0.933 | 0.908 | 1.778 |
|   |      | 1.5 | 0.113 | 2.534 | 1.899 | 8.319 |
|   | 0.5  | 0.5 | 0.059 | −0.041 | 0.031 | 0.033 |
|   |      | 1   | 0.083 | 1.002 | 1.093 | 2.098 |
|   |      | 1.5 | 0.111 | 2.542 | 1.924 | 8.388 |

**Table 2.** Monte Carlo simulations for $n = 400$, with optimal bandwidth given by (3.10)

| Model | $a_0$ | $c$ | Mean($h$) | 10 Bias | 100 Var | 100 MSE |
|---|---|---|---|---|---|---|
| 1 | 0.25 | 0.5 | 0.033 | 0.019 | 0.029 | 0.029 |
|   |      | 1   | 0.057 | 0.799 | 0.365 | 1.003 |
|   |      | 1.5 | 0.087 | 2.175 | 0.838 | 5.566 |
|   | 0.5  | 0.5 | 0.053 | −0.208 | 0.053 | 0.097 |
|   |      | 1   | 0.097 | 0.233 | 0.429 | 0.483 |
|   |      | 1.5 | 0.111 | 1.573 | 0.944 | 3.417 |
| 2 | 1    | 0.5 | 0.047 | −0.039 | 0.012 | 0.013 |
|   |      | 1   | 0.089 | 0.446 | 0.373 | 0.572 |
|   |      | 1.5 | 0.108 | 1.657 | 0.927 | 3.674 |
|   | 2    | 0.5 | 0.036 | 0.018 | 0.017 | 0.017 |
|   |      | 1   | 0.075 | 0.561 | 0.366 | 0.680 |
|   |      | 1.5 | 0.105 | 1.723 | 0.856 | 3.823 |
| 3 | 0.25 | 0.5 | 0.047 | −0.029 | 0.011 | 0.012 |
|   |      | 1   | 0.084 | 0.485 | 0.360 | 0.595 |
|   |      | 1.5 | 0.107 | 1.691 | 0.921 | 3.782 |
|   | 0.50 | 0.5 | 0.045 | −0.087 | 0.027 | 0.034 |
|   |      | 1   | 0.078 | 0.478 | 0.384 | 0.613 |
|   |      | 1.5 | 0.106 | 1.674 | 0.909 | 3.712 |

In the second simulation study, the estimator $\hat{a}(x)$ is compared with

$$\tilde{a}^{\#}(x) = \sup\{Y_i : \|X_i - x\| \leq h\}.$$

This estimator has the advantage of being simpler than that developed in Section 2, although it can be shown to have the inferior convergence rate of $n^{-1/(p+c)}$, rather than $n^{-2/(p+2c)}$. As for the first study, 100 arbitrary samples of size $n = 200$ and $n = 400$ are generated. We set $x = 0.5$ and $c = 0.5, 1$ or $1.5$. In order to make a fair comparison of the two competitors, the bandwidth $h$ of either method is selected by minimising in a deterministic way (i.e., minimising over the 100 samples) the MSE of each estimator. When constructing $\hat{a}(x)$, the choice of $h_1$ and $K$ and the estimation of $\ddot{a}(x), b(x)$ and $c$ are done as for the first study. It is found that, in 29 out of the 36 cases treated in Tables 1 and 2, the mean square error of $\hat{a}(x)$ is less than that of $\tilde{a}^{\#}(x)$, and that the median value of the ratio of mean square errors equals 0.21.

### 4.2. Data analysis

We consider data on 123 American electric utility companies, studied by Christensen and Greene (1976), Greene (1990) and Hall and Simar (2002), among others. We focus here on the relation between $Y_i = -\log(C_i/P_i)$ and $X_i = \log(Q_i)$, where $C_i$ is the cost, $Q_i$ the output and $P_i$ the price of fuel for each company. We fit the model

$$Y_i = a(X_i) + \varepsilon_i,$$

where it is assumed that the conditional density of the errors $\varepsilon_i$ satisfies relation (2.2). The scatterplot of the data, together with the estimated frontier curve $\hat{a}(x)$, is shown in Figure 2. We restrict the region of estimation to $[4.6, 11.2]$, to avoid estimation in sparse areas of $X$. Both the estimation of $\tilde{a}(x)$ and $\hat{a}(x)$ is done using local linear smoothing. At each point of an equispaced grid of 34 values between 4.6 and 11.2 we estimate the bandwidth $h = h_1$ from formula (3.8), yielding values in the range from 0.77 to 1.31. The bandwidth for estimating $\ddot{a}(x)$, $b(x)$ and $c(x)$ is chosen as one-fifth of the total range, namely 1.32, whereas to estimate the design density we use kernel estimation based on the normal reference rule. The kernel used throughout is again the biquadratic kernel.

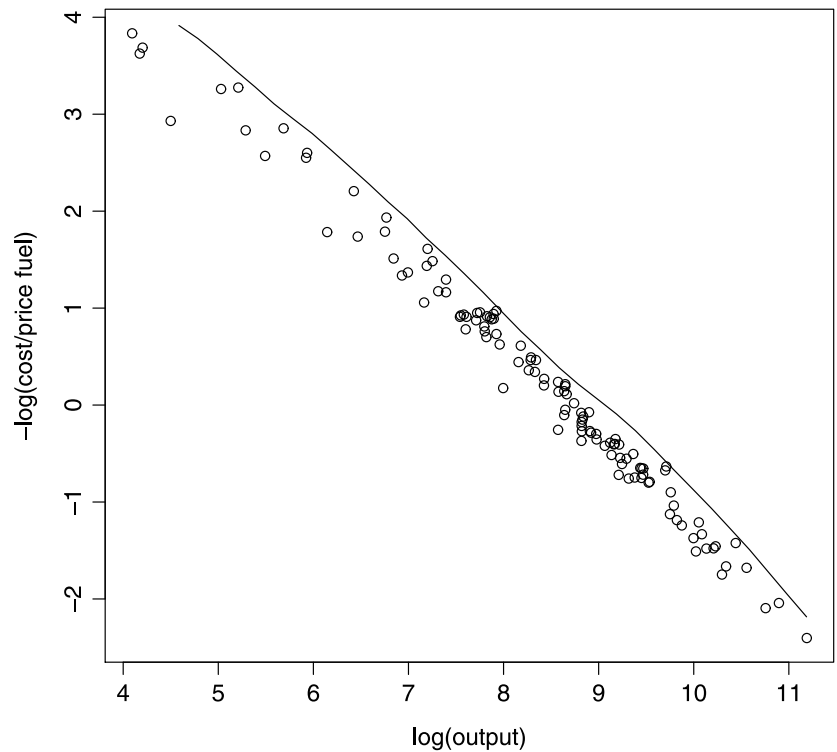

**Figure 2.** Scatterplot of the American Electric Utility Data. The observations are represented by circles, the solid curve is the estimated 'regression' (frontier) curve.

Figure 2 suggests that a linear model is appropriate for these data. However, it is particularly satisfying to reach that conclusion using a highly adaptive method that does not impose linearity, or even convexity, as a prior assumption.

## 5. Technical arguments

### 5.1. Proof of Theorem 1

To simplify notation we shall assume that $w_x = 1$ throughout; this can always be achieved via a change of scale. For brevity we shall deal only with the case where $x$ is an interior point of $\mathcal{R}$. Put $\gamma_\alpha(x) = a(x) - \alpha$, a scalar, and $\gamma_\beta(x) = h^{-1}\{\dot{a}(x) - \beta\}$, a $p$-vector. Let $\mathcal{I}(x, h)$ denote the set of indices $i$ such that $\|X_i - x\| \leq h$ and for $i \in \mathcal{I}(x, h)$ define $V_i = h^{-1}(X_i - x)$. In this notation,

$$Y_i - \alpha - \beta^{\mathrm{T}}(X_i - x) = \gamma_\alpha(x) + h^2\{\gamma_\beta(x)^{\mathrm{T}} V_i + \tfrac{1}{2} V_i^{\mathrm{T}} \ddot{a} V_i\} + h^2 R_i(x) + \varepsilon_i,$$

where the remainder, $R_i(x)$, has the property that

$$\begin{aligned}
\sup_{x \in \mathcal{R}} \sup_{i \in \mathcal{I}(x,h)} |R_i(x)| &\leq R(h) \\
&\equiv h^{-2} \sup_{x \in \mathcal{R}} \sup_{u: \|u\| \leq 1, x+hu \in \mathcal{R}} |a(x+hu) - a(x) - hu^{\mathrm{T}} \dot{a}(x) - \tfrac{1}{2} h^2 u^{\mathrm{T}} \ddot{a}(x) u|,
\end{aligned} \tag{5.1}$$

and $R(h) \to 0$ as $h \to 0$.

In particular, asking that $Y_i \geq \alpha + \beta^{\mathrm{T}}(X_i - x)$ for all indices $i \in \mathcal{I}(x, h)$ is equivalent to insisting that

$$\gamma_\alpha + \inf_{i \in \mathcal{I}(x,h)} \{h^2(\gamma_\beta^{\mathrm{T}} V_i + \tfrac{1}{2} V_i^{\mathrm{T}} \ddot{a} V_i) + h^2 R_i + \varepsilon_i\} \geq 0, \tag{5.2}$$

where we have dropped the argument from $\gamma_\alpha(x)$, $\gamma_\beta(x)$, $\ddot{a}(x)$ and $R_i(x)$. Let $\mathcal{S}_1(x, h)$ denote the set of pairs $(\gamma_\alpha, \gamma_\beta)$ such that (5.2) holds, and let $\tilde{\gamma}_1$ denote the infimum of $\gamma_\alpha$ over $(\gamma_\alpha, \gamma_\beta) \in \mathcal{S}_1(x, h)$. Then, $\tilde{a}(x) = a(x) - \tilde{\gamma}_1$.

It follows from this result and (5.1) that, defining

$$\tilde{\gamma}_2 = \tilde{\gamma}_2(x) = \sup_{\gamma_\beta} \inf_{i \in \mathcal{I}(x,h)} \{h^2(\gamma_\beta^{\mathrm{T}} V_i + \tfrac{1}{2} V_i^{\mathrm{T}} \ddot{a} V_i) + \varepsilon_i\} \tag{5.3}$$

and noting that, for any random variable $A$, essup $A$ is the infimum of constants $C$ for which $P(A \leq C) = 1$, we have:

$$h^{-2} \operatorname*{essup}_{x \in \mathcal{R}} \sup |\tilde{a}(x) - a(x) - \tilde{\gamma}_2| \to 0. \tag{5.4}$$

Defining $N = N(x, h) = \#\mathcal{I}(x, h)$, we may write $\tilde{\gamma}_2$ equivalently as

$$\tilde{\gamma}_2 = (nh^p)^{-1/c(x)} \sup_{\gamma_\beta} \inf_{1 \leq i \leq N} \{\rho_1(\gamma_\beta^{\mathrm{T}} V_{(i)} + \tfrac{1}{2} V_{(i)}^{\mathrm{T}} \ddot{a} V_{(i)}) + \xi_{(i)}\},$$

where $\rho_1 = (nh^p)^{1/c(x)} h^2$, $\xi_{(1)} < \xi_{(2)} < \cdots$ are the ordered values of $(nh^p)^{1/c(x)} \varepsilon_i$ for $i \in \mathcal{I}$ and $V_{(1)}, V_{(2)}, \ldots$ denote the concomitant values of $V_1, V_2, \ldots$. The factor $(nh^p)^{1/c(x)}$ here reflects the fact that the expected number of values $X_i$ that lie within $h$ of $x$ is asymptotic to a constant multiple of $nh^p$. The power, $1/c(x)$, is appropriate because it is the power applied to sample size to describe the scale of the largest value when the sample is drawn from a Pareto-type distribution with exponent $c(x)$.

For each $r \geq 1$,

the limiting joint distribution of $\xi_{(1)}, \ldots, \xi_{(r)}$ and $V_{(1)}, \ldots, V_{(r)}$ is the distribution of $b(x)^{-1/c(x)}(Z_1, \ldots, Z_r)$ and $U_1, \ldots, U_r$, where the sequence $Z_1, Z_2, \ldots$ is as defined at (3.5) and, independently of the $Z_j$'s, the $U_j$'s are uniformly distributed in the unit $p$-variate sphere. (5.5)

(See Hall (1978).) Moreover, with probability 1, for any interval $[a, b]$ where $0 < a < b < \infty$, the suprema over $\rho_1 \in [a, b]$ of the values of $\|\gamma\|$ and $i \geq 1$ at which the extremum

$$\sup_{\gamma \in \mathbb{R}^p} \inf_{1 \leq i < \infty} [\rho_1\{\gamma^{\mathrm{T}} U_i + \tfrac{1}{2} U_i^{\mathrm{T}} \ddot{a}(x) U_i\} + b(x)^{-1/c(x)} Z_i]$$

is achieved are finite and, using (5.5), the same can be proved of the extremum

$$\sup_{\gamma \in \mathbb{R}^p} \inf_{1 \leq i < \infty} \{\rho_1(\gamma^{\mathrm{T}} V_{(i)} + \tfrac{1}{2} V_{(i)}^{\mathrm{T}} \ddot{a} V_{(i)}) + \xi_{(i)}\}.$$

These properties and (5.5) imply that, if $\rho_1 \to \rho \in (0, \infty)$ as $n \to \infty$,

$$(nh^p)^{1/c(x)} \tilde{\gamma}_2 \to \sup_{\beta \in \mathbb{R}^p} \inf_{1 \leq i < \infty} [\rho\{\beta^{\mathrm{T}} U_i + \tfrac{1}{2} U_i^{\mathrm{T}} \ddot{a}(x) U_i\} + b(x)^{-1/c(x)} Z_i]$$

in distribution. The part of Theorem 1 pertaining to $\rho_1 \to \rho \in (0, \infty)$ follows from this property and (5.4).

If $\rho_1 \to 0$ then, since $\xi_{(1)} \to b(x)^{-1/c(x)} Z_1$ in distribution, we have $(nh^p)^{1/c(x)} \tilde{\gamma}_2 \to b(x)^{-1/c(x)} Z_1$ in distribution. And if $\rho_1 \to \infty$ then

$$h^{-2} \tilde{\gamma}_2 \to \sup_{\beta \in \mathbb{R}^p} \inf_{1 \leq i < \infty} \{\beta^{\mathrm{T}} U_i + \tfrac{1}{2} U_i^{\mathrm{T}} \ddot{a}(x) U_i\} = \sup_{\beta \in \mathbb{R}^p} \inf_{\|u\| \leq 1} \{\beta^{\mathrm{T}} u + \tfrac{1}{2} u^{\mathrm{T}} \ddot{a}(x) u\},$$

a constant. Parts (a) and (b) of Theorem 1 are consequences of these properties.

To establish convergence of second moments it suffices, in view of (5.4), to prove that for some $\eta_1 > 0$,

> there exist random variables $A_1$ and $A_2$ such that $A_1 \leq (nh^p)^{1/c(x)}\tilde{\gamma}_2 \leq A_2$ with probability 1, and $E(|A_j|^{2+\eta_1})$ is uniformly bounded for $j = 1, 2$. (5.6)

A proof of (5.6) is given in a longer version of this paper, obtainable from the authors.

## 5.2. Proof of Theorem 2

(Recall that we assume that $w_x = 1$.) We shall work with the definition (3.7) of $\check{a}(x)$. Defining $\tilde{\gamma}_2(x)$ as at (5.3), and noting (5.4), we have:

$$\check{a}(x) = \int a(x + h_1 u) K(u) \, du + \int \tilde{\gamma}_2(x + h_1 u) K(u) \, du + o_p(h^2). \quad (5.7)$$

The first integral on the right-hand side, $I_1(x)$, equals $a(x) + h^2 g(x) + o(h^2)$, where $g(x) = \frac{1}{2} t^2 \kappa (\nabla^2 a)(x)$, whence it follows that $(nh^p)^{1/c(x)}\{I_1(x) - a(x)\} \to \rho g(x)$. The stochastic process $S(u) = (nh^p)^{1/c(x)} \tilde{\gamma}_2(x + h_1 u)$ converges weakly to $Q_2(u) \equiv Q_2(u|x)$ (see below), whence it follows that the second integral, $I_2(x)$ on the right-hand side of (5.7), satisfies $(nh^p)^{1/c(x)} I_2(x) \to \int Q_2(u) K(u) \, du$.

To appreciate why the finite-dimensional distributions of $S$ converge to those of $Q_2$, consider the marked point process in $\mathbb{R}^d$, where the $i$th point is $V_i = h^{-p}(X_i - x)$ and the associated mark is $\zeta_i = \{nh^p(t+1)^p\}^{1/c(x)} \varepsilon_i$. Only the marked points that lie in the disc of radius $t+1$, centred at $x$, contribute to $\check{a}(x)$, and so we confine attention to those. Define $\zeta_{(1)} < \zeta_{(2)} < \cdots$ to be the ordered values of $\zeta_1 < \zeta_2 < \cdots$, and let $V_{(1)}, V_{(2)}, \ldots$ be the concomitant values of $V_1, V_2, \ldots$. In this new notation, (5.5) continues to hold. From that result it follows, using the argument in the paragraph containing (5.5), that for each finite set $u_1, \ldots, u_k$ in the sphere of radius $t+1$, centred at $x$, the joint distribution of $S(u_1), \ldots, S(u_k)$ converges to that of $Q_2(u_1), \ldots, Q_2(u_k)$. Tightness of the stochastic process $S$ can be proved using the fact that, defining

$$D(u, j_0) = \sup_{\gamma_\beta} \inf_{1 \leq j \leq j_0} \{ \rho_1 (\gamma_\beta^{\mathrm{T}} V_{(i_j(u))} + \tfrac{1}{2} V_{(i_j(u))}^{\mathrm{T}} \ddot{a} V_{(i_j(u))}) + (1+p)^{-p/c(x)} \zeta_{(i_j(u))} \},$$

where the ordering $j_1(u), j_2(u), \ldots$ is such that $V_{(i_1(u))} < V_{(i_2(u))} < \cdots$ among all indices $i(u)$ such that $\|V_{(i(u))} - u\| \leq 1$, the process $D(\cdot, j_0)$ decreases with increasing $j_0$.

## 5.3. Proof of Theorem 3

Derivation of (3.13) is similar to that of the last part of Theorem 1, and so will not be given here. We shall outline proofs of (3.14) and (3.15).

In the case of (3.14), take $a(x) = \delta^2 \psi(x/\delta)$ where $\delta = n^{-1/(p+2c)}$ and $\psi$ is a spherically symmetric function supported on $\mathcal{V}(0, \tfrac{1}{2})$ with bounded derivatives of first and second orders, all of them dominated by $\frac{1}{2}C$. Then, $a \in \mathcal{A}$. Consider the problem of discriminating between the models (a) $Y_i = \varepsilon_i$ and (b) $Y_i = a(X_i) + \varepsilon_i$ using only the data $(X_1, Y_1), \ldots, (X_n, Y_n)$. The likelihood-ratio approach, which in view of the Neyman–Pearson lemma is optimal, is to decide in favour of model (b) if and only if the ratio

$$L = \prod_{i=1}^{n} [f\{Y_i - a(X_i)\}/f(Y_i)]$$

exceeds an appropriate critical point. Here, $f$ is the density at (3.11). If $Y_i = Ia(X_i) + \varepsilon_i$, where $I = 0$ or 1 in cases (a) or (b), respectively, then

$$\log L = (c-1) \sum_{i=1}^{n} \log\{1 - Y_i^{-1} a(X_i)\} + \sum_{i=1}^{n} a(X_i).$$

Hence, the likelihood-ratio rule involves deciding in favour of (b) if and only if the sum $\ell = \sum_i \log\{1 - Y_i^{-1} a(X_i)\}$ exceeds a critical point.

Asymptotically correct discrimination is readily seen to be impossible if $\nu \equiv n\delta^p$ is bounded; this quantity is of the same order as the number of pairs $(X_i, Y_i)$ for $1 \leq i \leq n$ that contain information about $a$. The theorem will follow if we show that, when $\nu \equiv n\delta^p \to \infty$ but $\delta = o(n^{-2/(p+2c)})$ along a subsequence, the probability of correct discrimination using the likelihood-ratio rule when cases (a) and (b) above both have prior probability $\frac{1}{2}$ converges to $\frac{1}{2}$; it is assumed that all calculations are done for the subsequence.

We may Taylor-expand $\ell$, showing that $\ell/\delta^2 = \ell_1 + (I - \frac{1}{2})\ell_2 + \ell_3$, where $\ell_1 = -\sum_i \varepsilon_i^{-1} \psi_i$, $\ell_2 = \delta^2 \sum_i \varepsilon_i^{-2} \psi_i^2$, $\psi_i = \psi(X_i/\delta)$ and, when $\nu \to \infty$ and $\delta = o(n^{-1/(p+2c)})$, the remainder, $\ell_3$, equals $o_p(|\ell_1| + |\ell_2|)$. Using the fact that $0 < c < 2$ it may be proved that $\nu^{-1/c} \ell_1$ has a limiting, symmetric, non-degenerate stable distribution with exponent $c$, and $\delta^{-2} \nu^{-2/c} \ell_2$ has a limiting, positive, non-degenerate stable law with exponent $c/2$. Therefore, if $\delta = o(n^{-1/(p+2c)})$ then $\ell_2 = o_p(\ell_1)$, from which it follows that the probability of correct classification using the likelihood-ratio rule converges to $\frac{1}{2}$.

To obtain (3.15), let $\mathcal{W}$ denote the cube of diameter 2 inscribed within $\mathcal{V}(0,1)$, with its sides parallel to the coordinate axes. Place into $\mathcal{W}$ a rectangular grid of points, $x_1, \ldots, x_N$ with nearest neighbours exactly $\delta$ apart and no point distant less than $\frac{1}{2}\delta$ from the boundary of $\mathcal{V}(0,1)$. We may take $N \sim const.\delta^{-p}$ as $\delta \to 0$. Define $a_I(x) = \delta^2 \sum_i I_i \psi\{(x - x_i)/\delta\}$, where $I = (I_1, \ldots, I_N)$ is a vector of 0's and 1's. Then $a_I \in \mathcal{A}$ for each choice of $I$. Since $\psi$ vanishes outside radius $\frac{1}{2}$ from the origin, for each $x$ no more than one term in this series is non-zero. Treating the problem of estimating $a_I$ on $\mathcal{R}_h$ as one of discriminating between $I_i = 0$ and $I_i = 1$ for each $i$ such that the sphere of radius $\frac{1}{2}\delta$ centred at $x_i$ intersects $\mathcal{R}_h$ and arguing as in the proof of (3.14) we may derive (3.15).

## 6. Conclusion

We have shown that an alternative "regression" problem, where errors are "positioned" at their end-points, leads to estimators with properties quite different from their coun-terparts in conventional problems. In particular, if the error density is bounded away from zero then relative fast convergence rates are possible, even if the regression mean is known only up to smoothness conditions. Results of this type can be compared with their counterparts in $L_1$ or $L_2$ regression, where convergence rates using either method can be faster depending primarily on properties of the error distribution. Particularly if the main object were to obtain an idea of the shape of the regression mean, rather than for formal prediction, it would be appropriate to exploit these dissimilarities and construct a function estimator that enjoyed good convergence rates. Potential future problems of interest include developing adaptive methods for choosing among different regression methods, suitable for different error types, so as to ensure good empirical performance.

## Acknowledgements

We are grateful to two reviewers for helpful comments. The research of Van Keilegom was supported by IAP research network grant nr. P5/24 of the Belgian government (Belgian Science Policy). The authors would like to thank Léopold Simar for helpful discussions and for providing the data.

*Nonparametric "regression" when errors are positioned at end-points* 635